\documentclass[a4paper]{article}

\usepackage{amsfonts, latexsym, amssymb, amsthm, amsmath, mathrsfs, esint}

\newcommand{\R}{\mathbb{R}}
\newcommand{\N}{\mathbb{N}}
\newcommand{\Z}{\mathbb{Z}}

\newcommand{\Ha}{\mathcal{H}}
\newcommand{\F}{\mathcal{F}}
\newcommand{\E}{\mathcal{E}}
\newcommand{\Null}{\mathcal{N}}
\newcommand{\bfnu}{\boldsymbol{\nu}}
\newcommand{\bfx}{\boldsymbol{x}}
\newcommand{\bfXi}{\boldsymbol{\Xi}}
\newcommand{\U}{\boldsymbol{U}}
\newcommand{\ug}{\underline{g}}
\newcommand{\og}{\overline{g}}

\newcommand{\BV}{\mathrm{BV}}

\newcommand{\loc}{\mathrm{loc}}

\DeclareMathOperator{\dist}{dist}

\newcommand{\blank}{{\mkern 2mu\cdot\mkern 2mu}}
\DeclareMathOperator{\graph}{graph}
\newcommand{\Jac}{J}

\newcommand{\dd}[2]{\frac{\partial #1}{\partial #2}}
\newcommand{\set}[2]{\left\{ #1 \colon #2 \right\}}

\newcommand{\restr}{\mathchoice
{\kern2pt\mbox{\vrule width 0.08ex height1.5ex depth0ex\kern-0.08ex\vrule width 1.5ex height.08ex depth0ex}\kern2pt}
{\kern2pt\mbox{\vrule width 0.08ex height1.5ex depth0ex\kern-0.08ex\vrule width 1.5ex height.08ex depth0ex}\kern2pt}
{\kern1.5pt\mbox{\vrule width 0.06ex height1.1ex depth0ex\kern-0.06ex\vrule width 1.1ex height.06ex depth0ex}\kern1.5pt}
{\kern1pt\mbox{\vrule width 0.04ex height0.75ex depth0ex\kern-0.04ex\vrule width 0.75ex height.04ex depth0ex}\kern1pt}
}

\newtheorem{theorem}{Theorem}
\newtheorem{lemma}[theorem]{Lemma}
\newtheorem{proposition}[theorem]{Proposition}
\newtheorem{definition}[theorem]{Definition}

\theoremstyle{remark}
{}

\begin{document}

\title{Regularity of gradient vector fields giving rise to finite Caccioppoli partitions}

\author{Roger Moser\footnote{Department of Mathematical Sciences,
University of Bath, Bath BA2 7AY, UK.
E-mail: r.moser@bath.ac.uk}}

\maketitle

\begin{abstract}
For a finite set $A \subseteq \R^n$, consider a function $u \in \BV_\loc^2(\R^n)$ such that
$\nabla u \in A$ almost everywhere. If $A$ is convex independent, then it follows that $u$ is piecewise
affine away from a closed, countably $\Ha^{n - 1}$-rectifiable set. If $A$ is affinely independent,
then $u$ is piecewise affine away from a closed $\Ha^{n - 1}$-null set.

\end{abstract}

\section{Introduction}

For $n \in \N$, consider a finite set $A \subseteq \R^n$.
We study continuous functions $u \colon \R^n \to \R$ such that
the weak gradient $\nabla u$ satisfies $\nabla u \in \BV_\loc(\R^n; \R^n)$ and
$\nabla u(x) \in A$ for almost every $x \in \R^n$.
This means that whenever $\Omega \subseteq \R^n$ is open and bounded,
the sets $\set{x \in \Omega}{\nabla u(x) = a}$, for $a \in A$,
form a Caccioppoli partition of $\Omega$ as discussed, e.g., by Ambrosio, Fusco, and Pallara
\cite[Section 4.4]{Ambrosio-Fusco-Pallara:00}. The theory of Caccioppoli partitions therefore applies and gives some information
on the structure of $\nabla u$ and of $u$. The fact that we are dealing with
a gradient, however, gives rise to a better theory, especially under additional
assumptions on the geometry of $A$. We work with the following notions in this paper.

\begin{definition}
A set $A \subset \R^n$ is called \emph{convex independent} if any $a \in A$ does not
belong to the convex hull of $A \setminus \{a\}$. It is called \emph{affinely independent}
if any $a \in A$ does not belong to the affine span of $A \setminus \{a\}$.
\end{definition}

If either of these conditions is satisfied, then we can prove statements on
the regularity of $u$ that finite Caccioppoli partitions do not share in general.
In fact, we will see that $u$
is locally piecewise affine away from a closed, countably $\Ha^{n - 1}$-rectifiable set (if $A$ is convex independent)
or away from a closed $\Ha^{n - 1}$-null set (if $A$ is affinely independent).

In order to make this more precise, we introduce some notation.
Given $r > 0$ and $x \in \R^n$, we write
$B_r(x)$ for the open ball of radius $r$ centred at $x$.
Given $a \in \R^n$, the function $\lambda_a \colon \R^n \to \R$ is defined by $\lambda_a(x) = a \cdot x$
for $x \in \R^n$. Given two functions $v, w \colon \R^n \to \R$, we write $v \wedge w$ and $v \vee w$,
respectively, for the functions with $(v \wedge w)(x) = \min\{v(x), w(x)\}$ and $(v \vee w)(x) = \max\{v(x), w(x)\}$
for $x \in \R^n$.

\begin{definition} \label{def:regular}
Given a function $u \colon \R^n \to \R$, the \emph{regular set} of $u$, denoted by $\mathcal{R}(u)$, consists
of all $x \in \R^n$ such that there exist $a, b \in \R^n$, $c \in \R$, and $r > 0$ with
$u = \lambda_a \wedge \lambda_b + c$ in $B_r(x)$ or $u = \lambda_a \vee \lambda_b + c$ in $B_r(x)$.
The \emph{singular set} of $u$ is its complement $\mathcal{S}(u) = \R^n \setminus \mathcal{R}(u)$.
\end{definition}

The condition for $\mathcal{R}(u)$ allows the possibility that $a = b$, in which case $u$ is affine near $x$.
If $a \neq b$, then it is still piecewise affine near $x$.
Obviously $\mathcal{R}(u)$ is an open set and $\mathcal{S}(u)$ is closed.

It would be reasonable to include functions consisting of more than two affine pieces in the definition of
$\mathcal{R}(u)$, for example $(\lambda_{a_1} \wedge \lambda_{a_2}) \vee \lambda_{a_3} + c$ for $a_1, a_2, a_3 \in \R^n$ and $c \in \R$.
For the results of this paper, however, this would make no difference, therefore we choose the
simpler definition.

For $s \ge 0$, we denote the $s$-dimensional Hausdorff measure in $\R^n$ by $\Ha^s$. 
The notation $\BV_\loc^2(\R^n)$ is used for the space of functions
with weak gradient in $\BV_\loc(\R^n; \R^n)$. Thus the hypotheses of the following
theorems are identical to the assumptions at the beginning of the introduction.

\begin{theorem} \label{thm:convex-independent}
Suppose that $A$ is a finite, convex independent set. Let $u \in \BV_\loc^2(\R^n)$ with $\nabla u(x) \in A$
for almost every $x \in \R^n$. Then $\mathcal{S}(u)$ is countably $\Ha^{n - 1}$-rectifiable.
\end{theorem}

\begin{theorem} \label{thm:affinely-independent}
Suppose that $A$ is a finite, affinely independent set. Let $u \in \BV_\loc^2(\R^n)$ with $\nabla u(x) \in A$
for almost every $x \in \R^n$. Then $\Ha^{n - 1}(\mathcal{S}(u)) = 0$.
\end{theorem}

For $n = 2$, Theorem \ref{thm:affinely-independent} was proved in a previous paper \cite{Moser:17}.
For higher dimensions, the result is new. Theorem \ref{thm:convex-independent} is new even for
$n = 2$. For $n = 1$, both statements are easy to prove.

The results are optimal in terms of the Hausdorff measures involved. Furthermore, the assumption of
convex/affine independence is necessary. Indeed, there are examples
of finite sets $A \subseteq \R^2$ and functions $u \in \BV_\loc^2(\R^2)$ with $\nabla u(x) \in A$ almost everywhere such that
\begin{itemize}
\item $\Ha^2(\mathcal{S}(u)) > 0$; or
\item $\Ha^1(\mathcal{S}(u)) > 0$ and $A$ is convex independent; or
\item $\Ha^s(\mathcal{S}(u)) = \infty$ for any $s < 1$ and $A$ is affinely independent.
\end{itemize}
All of these can be found in the author's previous paper \cite{Moser:17}.

Apart from being of obvious geometric interest, functions as described above appear
in problems from materials science. They naturally arise as limits in $\Gamma$-convergence theories in the spirit of Modica and Mortola \cite{Modica-Mortola:77.1, Modica-Mortola:77.2}
for quantities such as
\begin{equation} \label{eqn:Gamma-limit}
\int_\Omega \left(\epsilon |\nabla^2 u|^2 + \frac{W(\nabla u)}{\epsilon}\right) \, dx,
\end{equation}
where $\Omega \subseteq \R^n$ is an open set and $W \colon \R^n \to [0, \infty)$
is a function with $A = W^{-1}(\{0\})$. Functionals of this sort appear in certain models for the surface energy of nanocrystals
\cite{Stewart-Goldenfeld:92, Liu-Metiu:93,  Watson-Norris:06}. For $\Omega \subseteq \R^2$, functions $u \in \BV^2(\Omega)$ with
$\nabla u \in \{(\pm 1, 0), (0, \pm 1)\}$ have also been used by
Cicalese, Forster, and Orlando \cite{Cicalese-Forster-Orlando:19} for a
different sort of $\Gamma$-limit arising from a model for frustrated spin systems.

Functionals similar to \eqref{eqn:Gamma-limit}, but for maps $u \colon \Omega \to \R^n$, also appear in certain models for phase
transitions in elastic materials (see, e.g., the seminal paper of Ball and James \cite{Ball-James:87}
or the introduction into the theory by M\"uller \cite{Mueller:96}).
In this context, due to the frame indifference of the underlying models, the set
$W^{-1}(\{0\})$ is typically not finite. Sometimes, however, the frame indifference is disregarded
(as in the paper by Conti, Fonseca, and Leoni \cite{Conti-Fonseca-Leoni:02}), or the theory gives
a limit with $\nabla u \in \BV(\Omega; A)$ for a finite set $A \subseteq \R^{n \times n}$ anyway
(such as in recent results of Davoli and Friedrich \cite{Davoli-Friedrich:20, Davoli-Friedrich:20.2}).
In such a case, Theorem \ref{thm:convex-independent} and Theorem \ref{thm:affinely-independent}
are potentially useful, as they apply to the components (or other one-dimensional projections) of $u$.

In the proof of Theorem \ref{thm:affinely-independent}, we use some of the
tools from the author's previous paper \cite{Moser:17}. In particular, we will analyse the intersections of the graph of
$u$ with certain hyperplanes in $\R^{n + 1}$. We will see that these intersections
correspond to the graphs of functions with $(n - 1)$-dimensional domains and with properties similar to $u$.
The key ideas from the previous paper, however, are specific to $\R^2$, so we eventually use different
arguments. In this paper, we use the theory of $\BV_\loc(\R^n; \R^n)$ to a much greater extent.
The central argument will consider approximate jump points of $\nabla u$. Near such a point,
we know that $u$ is close to a piecewise affine function
in a measure theoretic sense by definition. We then use an induction
argument (with induction over $n$) to show that $u$ is in fact piecewise affine near $\Ha^{n - 1}$-almost every
approximate jump point.

We also need to analyse points where $u$ has an approximate limit, and they are of interest for the proofs of
both Theorem \ref{thm:affinely-independent} and Theorem \ref{thm:convex-independent}. This part of the analysis
is significantly simpler and relies on
the fact that for any $a \in A$, the function $v(x) = u(x) - a \cdot x$ has some monotonicity
properties.

In the rest of the paper, we study a fixed function $u \in \BV_\loc^2(\R^n)$ with $\nabla u(x) \in A$ for almost every $x \in \R^n$.
Since we are interested only in the local properties of $u$, we may assume that it is also bounded.
(Otherwise we can modify it outside of a bounded set with the construction described in \cite[Section 6]{Moser:17}.)
We define the function $\U \colon \R^n \to \R^{n + 1}$ by
\[
\U(x) = \begin{pmatrix} x \\ u(x) \end{pmatrix}, \quad x \in \R^n.
\]
We use the notation $\graph(u) = \U(\R^n)$ for the graph of $u$.

As we sometimes work with points in $\R^{n + 1}$ (especially points on $\graph(u)$) and their projections onto
$\R^n$ simultaneously, we use the following notation. A generic point in $\R^{n + 1}$ is denoted by $\bfx = (x_1, \dotsc, x_{n + 1})^T$,
and then we write $x = (x_1, \dotsc, x_n)^T$. Thus $\bfx = (\begin{smallmatrix} x \\ x_{n + 1} \end{smallmatrix})$.
We think of elements of $\R^n$ and of $\R^{n + 1}$
as column vectors, and this is sometimes important, as we use them as columns in certain matrices.

As our function satisfies in particular the condition $\nabla u \in \BV_\loc(\R^n; \R^n)$, the theory
of this space will of course be helpful. In this context, we mostly follow the notation and terminology
of Ambrosio, Fusco, and Pallara \cite{Ambrosio-Fusco-Pallara:00}. We also use several of the results found in this book.

\section{Approximate faces and edges of the graph} \label{sct:faces-edges}

In this section, we decompose $\R^n$ into three sets $\F$, $\E$, and $\Null$.
These are defined such that we expect regularity in $\F$ under the assumptions of either of the main theorems, and
also in $\E$ under the assumptions of Theorem \ref{thm:affinely-independent}. The third set, $\Null$, will be an
$\Ha^{n - 1}$-null set. The sets $\F$ and $\E$ characterised, up to $\Ha^{n - 1}$-null sets, by the condition that
$\nabla u$ has an approximate limit or an approximate jump, respectively. Since much of our analysis
examines $\graph(u)$, it is also convenient to think of $\F$ as the set of points where the graph behaves approximately
like the ($n$-dimensional) faces of a polyhedral surface, whereas $\E$ corresponds to approximate ($(n - 1)$-dimensional) edges.

First we define the set $\F' \subseteq \R^n$, comprising all points $x \in \R^n$ such that there exists
$a \in \R^n$ satisfying
\[
\lim_{r \searrow 0} \fint_{B_r(x)} |\nabla u - a| \, d\Ha^n = 0.
\]
In other words, this is the set of all points where $\nabla u$ has an approximate limit $a$.
It is then clear that $a \in A$. The complement $\R^n \setminus \F$ is called the approximate discontinuity set of $\nabla u$.

Furthermore, let $\E'$ be the set of all $x \in \R^n$ such that there exist $a_-, a_+ \in \R^n$
with $a_- \neq a_+$ and there exists $\eta \in S^{n - 1}$ such that
\begin{equation} \label{eqn:jump+}
\lim_{r \searrow 0} \fint_{\set{\tilde{x} \in B_r(x)}{(\tilde{x} - x) \cdot \eta > 0}} |\nabla u - a_+| \, d\Ha^n = 0
\end{equation}
and
\begin{equation} \label{eqn:jump-}
\lim_{r \searrow 0} \fint_{\set{\tilde{x} \in B_r(x)}{(\tilde{x} - x) \cdot \eta < 0}} |\nabla u - a_-| \, d\Ha^n = 0.
\end{equation}
This is the approximate jump set of $\nabla u$.
Again, the points $a_-, a_+$ will always belong to $A$.

According to a result by Federer and Vol'pert (which can be found in the book by Ambrosio, Fusco, and Pallara
\cite[Theorem 3.78]{Ambrosio-Fusco-Pallara:00}), there exists an $\Ha^{n - 1}$-null
set $\Null' \subseteq \R^n$ such that
\[
\R^n = \F' \cup \E' \cup \Null'.
\]
Furthermore, the set $\E'$ is countably $\Ha^{n - 1}$-rectifiable.

Given $x \in \R^n$ and $\rho > 0$, we define the function $u_{x, \rho} \colon \R^n \to \R$ with
\[
u_{x, \rho}(\tilde{x}) = \frac{1}{\rho}\left(u(x + \rho\tilde{x}) - u(x)\right)
\]
for $\tilde{x} \in \R^n$. For $x$ fixed, the family of functions $(u_{x, \rho})_{\rho > 0}$ is clearly bounded
in $C^{0, 1}(K)$ for any compact set $K \subseteq \R^n$. Therefore, the theorem of Arzel\`a--Ascoli implies that
there exists a sequence $\rho_k \searrow 0$ such that $u_{x, \rho_k}$ converges locally uniformly.
If we have in fact a limit for $\rho \searrow 0$, then we write
\[
T_x u = \lim_{\rho \searrow 0} u_{x, \rho}
\]
and call this limit the \emph{tangent function} of $u$ at $x$.

If $x \in \F'$ and $a \in A$ is the approximate limit of $\nabla u$ at $x$, then for any sequence $\rho_k \searrow 0$,
the limit of $u_{x, \rho_k}$ can only be $\lambda_a$. Hence in this case, there
exists a tangent function $T_x u$, which is exactly this function.
Similarly, if $x \in \E'$, then $T_x u$ exists and
\[
T_x u(\tilde{x}) =  \begin{cases}
\lambda_{a_-}(\tilde{x}) & \text{if $\tilde{x} \cdot \eta < 0$}, \\
\lambda_{a_+}(\tilde{x}) & \text{if $\tilde{x} \cdot \eta \ge 0$}.
\end{cases}
\]
Because $T_x u$ is a continuous function, this means that
\[
\eta = \pm \frac{a_+ - a_-}{|a_+ - a_-|}.
\]
Then we conclude that $T_x u = \lambda_{a_-} \wedge \lambda_{a_+}$ or $T_x u = \lambda_{a_-} \vee \lambda_{a_+}$,
depending on the sign.

If we consider the functions $a_-, a_+ \colon \E' \to A$ and $\eta \colon \E' \to S^{n - 1}$ such that
\eqref{eqn:jump+} and \eqref{eqn:jump-} are satisfied on $\E'$, then the previously used result
\cite[Theorem 3.78]{Ambrosio-Fusco-Pallara:00} also implies that
\[
D\nabla u \restr \E' = (a_+ - a_-) \otimes \eta \, \Ha^{n - 1} \restr \E'.
\]
Let $\gamma = \min\set{|a - b|}{a, b \in A}$. Then for any Borel set $\Omega \subseteq \R^n$, we conclude that
\[
|D\nabla u|(\Omega) \ge \gamma \Ha^{n - 1}(\E' \cap \Omega).
\]

Now define
\[
\F = \set{x \in \F'}{\lim_{\rho \searrow 0} \rho^{1 - n} |D\nabla u|(B_r(x)) = 0}.
\]
Then standard results \cite[Theorem 2.56 and Lemma 3.76]{Ambrosio-Fusco-Pallara:00} imply that 
$\Ha^{n - 1}(\F' \setminus \F) = 0$.

Recall the map $\U \colon \R^n \to \R^{n + 1}$ defined in the introduction.
Set $\F^* = \U(\F)$ and $\E^\dagger = \U(\E')$. Then $\E^\dagger$ is
a countably $\Ha^{n - 1}$-rectifiable subset of $\R^{n + 1}$. Hence at $\Ha^{n - 1}$-almost
every $\bfx \in \E^\dagger$, the measure $\Ha^{n - 1} \restr \E^\dagger$ has a
tangent measure \cite[Theorem 2.83]{Ambrosio-Fusco-Pallara:00} of the form $\Ha^{n - 1} \restr T_{\bfx} \E^\dagger$, where
$T_{\bfx} \E^\dagger$ is an $(n - 1)$-dimensional linear subspace of $\R^{n + 1}$
(the approximate tangent space of $\E^\dagger$ at $\bfx$). Let $\E^*$ be the set
of all $\bfx \in \E^\dagger$ where this is the case. Furthermore, let
$\E = \U^{-1}(\E^*)$. Then $\E' \setminus \E$ is an $\Ha^{n - 1}$-null set.

Thus if we define $\Null = \R^n \setminus (\F \cup \E)$, then
$\Null$ is an $\Ha^{n - 1}$-null set and we have the disjoint decomposition
\[
\R^n = \F \cup \E \cup \Null.
\]

\section{Proof of Theorem \ref{thm:convex-independent}}

In this section we prove our first main result, Theorem \ref{thm:convex-independent}.
The proof is based on the following proposition, which will also be useful for the proof of Theorem \ref{thm:affinely-independent} later on.

\begin{proposition} \label{prp:convex-independent}
Suppose that $A \subseteq \R^n$ is finite and convex independent. Let $u \in \BV_\loc^2(\R^n)$
be a function with $\nabla u(x) \in A$ for almost all $x \in \R^n$.
Then there exist $r > 0$ and $\epsilon > 0$ with the following property.
Suppose that there exists $a \in A$ such that
\begin{equation} \label{eqn:measure-nabla-u-neq-a}
\Ha^n(\set{x \in B_1(0)}{\nabla u(x) \neq a}) \le \epsilon
\end{equation}
and
\[
|D\nabla u|(B_1(0)) \le \epsilon.
\]
Then $\nabla u(x) = a$ for almost every $x \in B_r(0)$.
\end{proposition}

\begin{proof}
Because $A$ is convex independent, there exists $\omega \in S^{n - 1}$ such that
\[
a \cdot \omega < \min_{b \in A \setminus \{a\}} b \cdot \omega.
\]
As $A$ is finite, there also exists $\delta \in (0, 1)$ such that the inequality $a \cdot \xi \le \min_{b \in A \setminus \{a\}} (b \cdot \xi)$
holds even for $\xi$ in the cone
\[
C = \set{\xi \in \R^n}{\xi \cdot \omega \ge \delta |\xi|}.
\]

Consider the function $v \colon \R^n \to \R$ with $v(x) = u(x) - a \cdot x$ for $x \in \R^n$.
Then for any $\xi \in C$,
\[
\xi \cdot \nabla v(x) = \xi \cdot \nabla u(x) - a \cdot \xi \ge 0
\]
almost everywhere. Thus $v$ is monotone along lines parallel to $\xi$.
(This is true for every such line by the continuity of $v$.)
Furthermore, for almost every $x \in \R^n$, we find that either $\nabla u(x) = a$
or $\omega \cdot \nabla v(x) > 0$.

Suppose that $\nabla u = a$ does \emph{not} hold almost everywhere in $B_r(0)$. Then
there exist $x_-, x_+ \in B_r(0)$ with $v(x_-) < v(x_+)$. Define
\[
C_- = (x_- - C) \cap B_1(0) \quad \text{and} \quad C_+ = (x_+ + C) \cap B_1(0).
\]
Then for any $x' \in C_-$ and $x'' \in C_+$, we conclude that
\[
v(x') \le x(x_-) < v(x_+) \le v(x'').
\]

We now foliate a part of $B_1(0)$ by line segments parallel to $\omega$.
For $R \in (0, 1]$, let $Z_R = \set{x \in B_R(0)}{\omega \cdot x = 0}$. For every $z \in Z_R$,
consider the line segment
\[
L_z = \set{z + t\omega}{-\frac{1}{2} \le t \le \frac{1}{2}}.
\]
Provided that $r$ is chosen sufficiently small, we can find $R \in (0, 1]$ such that
\[
\set{z - \frac{\omega}{2}}{z \in Z_R} \subseteq C_- \quad \text{and} \quad \set{z + \frac{\omega}{2}}{z \in Z_R} \subseteq C_+.
\]
Hence for any $z \in Z_R$,
\[
v\left(z + \frac{\omega}{2}\right) - v\left(z - \frac{\omega}{2}\right) \ge v(x_+) - v(x_-) > 0.
\]
In particular, the restriction of $v$ to the line segment $L_z$ is not constant. For $z \in Z_R$,
define $\Lambda_z = \set{x \in L_z}{\nabla u(x) = a}$. Then it follows that $\Ha^1(\Lambda_z) < 1$
for $\Ha^{n - 1}$-almost all $z \in Z_R$.

On the other hand, because of \eqref{eqn:measure-nabla-u-neq-a}, we also know that
\[
\Ha^{n - 1}\left(\set{z \in Z_R}{\Ha^1(\Lambda_z) = 0}\right) \le \epsilon.
\]
Thus if we define $Z' = \set{z \in Z_R}{0 < \Ha^1(\Lambda_z) < 1}$, then
\[
\Ha^{n - 1}(Z') \ge \Ha^{n - 1}(Z_R) - \epsilon.
\]

Set $c = \min_{b \in A} |a - b|$.
For $\Ha^{n - 1}$-almost any $z \in Z'$, the function $t \mapsto \nabla u(z + t\omega)$ belongs
to $\BV\bigl((-\frac{1}{2}, \frac{1}{2}); \R^n\bigr)$ and its total variation is at least $c$. Hence
\cite[Theorem 3.103]{Ambrosio-Fusco-Pallara:00}
\[
|D\nabla u|(B_1(0)) \ge c \Ha^{n - 1}(Z') \ge c\bigl(\Ha^{n - 1}(Z_R) - \epsilon\bigr).
\]
If $\epsilon$ is sufficiently small, then this means in particular that $|D\nabla u|(B_1(0)) > \epsilon$.
Thus we have proved the contrapositive of Proposition \ref{prp:convex-independent}.
\end{proof} 

\begin{proof}[Proof of Theorem \ref{thm:convex-independent}]
We show that $\F \subseteq \mathcal{R}(u)$. To this end, fix $x \in \F$ and consider the rescaled
functions $u_{x, \rho}$ for $\rho > 0$. Since $x \in \F$, we know that $\nabla u_{x, \rho} \to a$
in $L^1(B_1(0))$ as $\rho \searrow 0$ for some $a \in A$. Furthermore, since
\[
|D\nabla u_{x, \rho}|(B_1(0)) = \rho^{1 - n} |D\nabla u|(B_\rho(x)) \to 0
\]
as $\rho \searrow 0$, the function $u_{x, \rho}$ satisfies the inequalities of Proposition \ref{prp:convex-independent}
for $\rho$ sufficiently small. Hence $\nabla u_{x, \rho}(\tilde{x}) = a$ for almost every $\tilde{x} \in B_r(0)$,
which implies that
\[
u(\tilde{x}) = u(x) + a \cdot (\tilde{x} - x)
\]
for all $\tilde{x} \in B_{\rho r}(x)$. Hence $x \in \mathcal{R}(u)$.

Theorem \ref{thm:convex-independent} now follows from the observations in Section \ref{sct:faces-edges}.
\end{proof}

\section{Specialising to a regular $n$-simplex} \label{sct:simplex}

The rest of the paper is devoted to the proof of Theorem \ref{thm:affinely-independent}.
Instead of considering \emph{any} affinely independent set $A$, we now assume that
$a_0, \dotsc, a_n \in \R^n$ are the corners of a regular $n$-simplex of side length $\sqrt{2n + 2}$ centred
at $0$, and that $A = \{a_0, \dotsc, a_n\}$.
We further assume that the matrix with columns $a_0 - a_1, \dotsc, a_0 - a_n$ has a
positive determinant. Theorem \ref{thm:affinely-independent}
can then be reduced to this situation by composing $u$ with an affine transformation.
The details are given on page \pageref{reduction-to-simplex} below.

As it is sometimes convenient to permute $a_0, \dotsc, a_n$ cyclically, we
regard $0, \dotsc, n$ as members of $\Z_{n + 1} = \Z / (n + 1)\Z$ in this context. Thus $a_{i + n + 1} = a_i$.

The condition that our simplex has side length $\sqrt{2n + 2}$ means that $|a_i| = \sqrt{n}$
for every $i \in \Z_{n + 1}$. Indeed, by the calculations of Parks and Wills \cite{Parks-Wills:02},
the dihedral angle of the regular $n$-simplex is $\arccos \frac{1}{n}$. As each $a_i$ is orthogonal to
one of the faces, this means that
$a_i \cdot a_j = - \frac{1}{n} |a_i| |a_j|$ for $i \neq j$, and therefore
$2n + 2 = |a_i - a_j|^2 = \frac{2n + 2}{n} |a_i| |a_j|$. From this we conclude that
\[
|a_i| = \sqrt{n}
\]
for $i \in \Z_{n + 1}$ and
\[
a_i \cdot a_j = -1
\]
for $i \neq j$.

For $i \in \Z_{n + 1}$, we now define the vector $\bfnu_i \in \R^{n + 1}$ by
\[
\bfnu_i = \frac{1}{\sqrt{n + 1}} \begin{pmatrix} -a_i \\ 1 \end{pmatrix}.
\]
Then
\[
|\bfnu_i|^2 = \frac{|a_i|^2 + 1}{n + 1} = 1,
\]
whereas for $i \neq j$,
\[
\bfnu_i \cdot \bfnu_j = \frac{a_i \cdot a_j + 1}{n + 1} = 0.
\]
Hence $(\bfnu_1, \dotsc, \bfnu_{n + 1})$ is an orthonormal basis of $\R^{n + 1}$. (This is the reason why we
choose $A$ as above.) Furthermore,
\[
\begin{split}
\det \begin{pmatrix} -a_1 & \cdots & -a_{n + 1} \\ 1 & \cdots & 1 \end{pmatrix} & = \det \begin{pmatrix} a_0 - a_1 & \cdots & a_0 - a_n & -a_0 \\  0 & \cdots & 0 & 1 \end{pmatrix} \\
& = \det \begin{pmatrix} a_0 - a_1 & \cdots & a_0 - a_n \end{pmatrix}.
\end{split}
\]
(In the first step, we have used the fact that $a_{n + 1} = a_0$ and subtracted the last column from each
of the other columns of the matrix.) Hence
the above assumption guarantees that the basis $(\bfnu_1, \dotsc, \bfnu_{n + 1})$ gives the standard orientation of
$\R^{n + 1}$.

We now use the notation $\lambda_i = \lambda_{a_i}$, recalling that this is the linear function with
$\lambda_i(x) = a_i \cdot x$ for $x \in \R^n$. For $i \in \Z_{n + 1}$, we set
\[
\F_i = \set{x \in \F}{T_x u = \lambda_i}.
\]
Thus we have the disjoint decomposition
\[
\F = \bigcup_{i \in \Z_{n + 1}} \F_i.
\]
Furthermore, we define $\F_i^* = \U(\F_i)$.

Of course $\U \colon \R^n \to \graph(u)$ is a bi-Lipschitz map.
Thus in order to understand $\F$, $\E$, or $\F_i$, it suffices to study $\F^*$, $\E^*$, or $\F_i^*$
and how $\U^{-1}$ transforms them. In particular, the following is true.

\begin{lemma}
For any Borel set $\Omega \subseteq \R^n$,
\[
\Ha^{n - 1}(\E^* \cap (\Omega \times \R)) = \sqrt{\frac{n + 1}{2}} \Ha^{n - 1}(\E \cap \Omega) = \frac{1}{2} |D\nabla u|(\Omega).
\]
\end{lemma}

\begin{proof}
We use the area formula \cite[Theorem 2.91]{Ambrosio-Fusco-Pallara:00}. Hence
we need to calculate the Jacobian of $\U$ restricted to the approximate tangent spaces of $\E$.

More precisely, since $\E$ is countably $\Ha^{n - 1}$-rectifiable, there exists an approximate tangent space
$T_x \E$ at $\Ha^{n - 1}$-almost every $x \in \E$. Because $\U$ is Lipschitz continuous, the tangential
derivative $d^\E \U(x)$ exists at $\Ha^{n - 1}$-almost every $x \in \E$ \cite[Theorem 2.90]{Ambrosio-Fusco-Pallara:00}.
We write $L^*$ for the adjoint of a linear operator $L$. Then
\[
\Jac_\E \U(x) = \sqrt{\det\bigl((d^\E \U(x))^* \circ d^\E \U(x)\bigr)}
\]
is the Jacobian of $\U$ at $x$ with respect to $T_x\E$. The area formula implies that
\[
\Ha^{n - 1}(\U(\E \cap \Omega)) = \int_{\E \cap \Omega} \Jac_\E \U(x) \, d\Ha^{n - 1}.
\]
Thus in order to prove the first identity, it suffices to show that
\[
\Jac_\E \U(x) = \sqrt{\frac{n + 1}{2}}
\]
for $\Ha^{n - 1}$-almost every $x \in \E$.

To this end, consider $x \in \E$. Note that $T_x\E = (a_i - a_j)^\perp$ for some $i, j \in \Z_{n + 1}$ with $i \neq j$
at $\Ha^{n - 1}$-almost every such point. For $\xi \in (a_i - a_j)^\perp$, we know that
\[
\frac{1}{\rho}(u(x + \rho \xi) - u(x)) = u_{x, \rho}(\xi) \to T_x u(\xi)
\]
as $\rho \searrow 0$. The convergence is in fact uniform on compact subsets of $(a_i - a_j)^\perp$.
Moreover, since $T_x u = \lambda_i \wedge \lambda_j$ or $T_x u = \lambda_j \vee \lambda_j$, its restriction
to $(a_i - a_j)^\perp$ is linear with $T_x u(\xi) = a_i \cdot \xi$.
Hence $d^\E u(x)$ exists, and so does $d^\E \U(x)$. We calculate
\[
d^\E \U(x) \xi = \begin{pmatrix} \xi \\ a_i \cdot \xi\end{pmatrix}.
\]

For simplicity, we assume that $i = n - 1$ and $j = n$. The space $(a_i - a_j)^\perp$ is spanned
by the vectors $a_0, \dotsc, a_{n - 2}$. Suppose that we choose an orthonormal basis $(\epsilon_0, \dotsc, \epsilon_{n - 2})$
of $T_x\E$. Let $L \colon T_x\E \to T_x \E$ denote the linear operator that maps $\epsilon_i$ to $a_i$ for $i = 0, \dotsc, n - 2$.
Then $d^\E \U(x) \circ L$ is represented by the matrix
\[
M_1 = \begin{pmatrix} a_0 & \cdots & a_{n - 2} \\ a_0 \cdot a_{n - 1} & \cdots & a_{n - 2} \cdot a_{n - 1} \end{pmatrix} = \begin{pmatrix} a_0 & \cdots & a_{n - 2} \\ -1 & \cdots & -1 \end{pmatrix}
\]
with respect to the above basis. Hence
\[
\Jac_\E \U(x) = \sqrt{\frac{\det(M_1^T M_1)}{\det(L^* \circ L)}}.
\]
We write $I_k$ for the identity $k \times k$-matrix. Then
\[
M_1^T M_1 = \begin{pmatrix}
a_0 \cdot a_0 + 1 & \cdots & a_0 \cdot a_{n - 2} + 1 \\
\vdots & \ddots & \vdots \\
a_{n - 2} \cdot a_0 + 1 & \cdots & a_{n - 2} \cdot a_{n - 2} + 1
\end{pmatrix} = (n + 1) I_{n - 1}
\]
and $\det(M_1^T M_1) = (n + 1)^{n - 1}$.

As $L$ maps an $(n - 1)$-cube of side length $1$ to the parallelepiped spanned by $a_0, \dotsc, a_{n - 2}$,
we know that $\det(L^* \circ L)$ is the $(n - 1)$-volume of the latter. Thus if $M_2$ is the
$n \times (n - 1)$-matrix with columns $a_0, \dotsc, a_{n - 2}$, then
\[
\det(L^* \circ L) = \det(M_2^T M_2).
\]
We further compute
\[
M_2^T M_2 = \begin{pmatrix}
a_0 \cdot a_0 & \cdots & a_0 \cdot a_{n - 2} \\
\vdots & \ddots & \vdots \\
a_{n - 2} \cdot a_0 & \cdots & a_{n - 2} \cdot a_{n - 2}
\end{pmatrix} = \begin{pmatrix}
n & -1 & \cdots & -1 \\
-1 & n & & \vdots \\
\vdots & & \ddots & -1 \\
-1 & \cdots & -1 & n
\end{pmatrix}.
\]
In order to calculate the determinant, we first subtract the first row of this matrix from each of the
other rows. We obtain
\[
\begin{split}
\det(M_2^T M_2) & = \det \begin{pmatrix}
n & -1 & \cdots & \cdots & -1 \\
-(n + 1) & n + 1 & 0 & \cdots & 0 \\
-(n + 1) & 0 & n + 1 & & \vdots \\
\vdots & \vdots & & \ddots & 0 \\
-(n + 1) & 0 & \cdots & 0 & n + 1
\end{pmatrix} \\
& = (n + 1)^{n - 2} \det \begin{pmatrix}
n & -1 & \cdots & \cdots & -1 \\
-1 & 1 & 0 & \cdots & 0 \\
-1 & 0 & 1 & & \vdots \\
\vdots & \vdots & & \ddots & 0 \\
-1 & 0 & \cdots & 0 & 1
\end{pmatrix}.
\end{split}
\]
In the last matrix, we now add to the first row the sum of all the other rows. Thus
\[
\det(M_2^T M_2) = (n + 1)^{n - 2} \det \begin{pmatrix}
2 & 0 & \cdots & \cdots & 0 \\
-1 & 1 & 0 & \cdots & 0 \\
-1 & 0 & 1 & & \vdots \\
\vdots & \vdots & & \ddots & 0 \\
-1 & 0 & \cdots & 0 & 1
\end{pmatrix} = 2(n + 1)^{n - 2}.
\]
Hence
\[
\Jac_\E \U(x) = \sqrt{\frac{\det(M_1^T M_1)}{\det(M_2^T M_2)}} = \sqrt{\frac{n + 1}{2}}.
\]

In order to prove the second identity, we recall that $|a_i - a_j| = \sqrt{2n + 2}$ for $i \neq j$.
Hence $|D\nabla u|(\Omega) = \sqrt{2n + 2} \Ha^{n - 1}(\E \cap \Omega) = 2\Ha^{n - 1}(\E^* \cap (\Omega \times \R))$.
\end{proof}

\section{Slicing the graph}

We still assume that $A$ consists of the corners of the regular $n$-simplex from Section \ref{sct:simplex}
and we assume that $u \in \BV_\loc^2(\R^n)$ is bounded and satisfies $\nabla u(x) \in A$ for almost
every $x \in \R^n$. In this section, we analyse the graph of $u$. In particular, we examine intersections
of $\graph(u)$ with hyperplanes perpendicular to one of the vectors $\bfnu_i$. We will see that almost all
such intersections can be represented as the graphs of functions in $\BV_\loc^2(P)$, where
\[
P = \set{y \in \R^n}{y_1 + \dotsb + y_n = 0},
\]
and with gradient taking one of $n$ different values almost everywhere. That is, we have a function with
properties similar to $u$, but with an $(n - 1)$-dimensional domain.
This observation will eventually make it possible to prove Theorem \ref{thm:affinely-independent} with the
help of an induction argument.

We use some tools from the author's previous paper \cite{Moser:17} in this section.
Given $i \in \Z_{n + 1}$, let $\Phi_i \colon \R^{n + 1} \to \R^{n + 1}$ be the linear map with
\[
\Phi_i(\bfx) = \begin{pmatrix} \bfnu_{i + 1} \cdot \bfx \\ \vdots \\ \bfnu_{i + n + 1} \cdot \bfx \end{pmatrix},
\]
so that $\Phi_i(\bfnu_{i + k})$ is the $k$-th standard basis vector in $\R^{n + 1}$. For $t \in \R$, let
\[
\Gamma_i(t) = \set{y \in \R^n}{\begin{pmatrix} y \\ t \end{pmatrix} \in \Phi_i(\graph(u))}.
\]
This corresponds to the intersection of $\graph(u)$ with a hyperplane orthogonal to $\bfnu_i$ after
rotation by $\Phi_i$, or in other words, a slice of $\graph(u)$.

We further define the functions
\[
\ug_i(y) = \sup\set{t \in \R}{u(t\nu_i + y_1 \nu_{i + 1} + \dotsb + y_n \nu_{i + n}) > \frac{t + y_1 + \dotsb + y_n}{\sqrt{n + 1}}}
\]
and
\[
\og_i(y) = \inf\set{t \in \R}{u(t\nu_i + y_1 \nu_{i + 1} + \dotsb + y_n \nu_{i + n}) < \frac{t + y_1 + \dotsb + y_n}{\sqrt{n + 1}}}.
\]
Note that for a fixed $y \in \R^n$, the set
\[
\set{t \in \R}{u(t\nu_i + y_1 \nu_{i + 1} + \dotsb + y_n \nu_{i + n}) = \frac{t + y_1 + \dotsb + y_n}{\sqrt{n + 1}}}
\]
corresponds to the intersection of $\graph(u)$ with a line parallel to $\nu_i$, so the functions $\ug_i$
and $\og_i$ tell us something about the geometry of $\graph(u)$ as well.

The following properties of $\ug_i$ and $\og_i$ have been proved elsewhere for $n = 2$ \cite[Lemma 16]{Moser:17}.
The proof carries over to higher dimensions as well.
We therefore do not repeat it here.

\begin{lemma} \label{lem:g}
For any $i \in \Z_{n + 1}$, the following statements hold true.
\begin{enumerate}
\item The function $\ug_i$ is lower semicontinuous and $\og_i$ is upper semicontinuous.
\item The identity $\ug_i = \og_i$ holds almost everywhere in $\R^n$.
\item For any $y \in \R^n$, the inequality $\ug_i(y) \le \og_i(y)$ holds true and
\[
\{y\} \times [\ug_i(y), \og_i(y)] \subseteq \Phi_i(\graph(u)).
\]
\item \label{item:g.iv} Let $t \in \R$ and $y \in \R^n$. Then $y \in \Gamma_i(t)$ if, and only if, $\ug_i(y) \le t \le \og_i(y)$.
\item \label{item:g.v} For all $y \in \R^n$ and all $\zeta \in (0, \infty)^n$, the inequality $\og_i(y + \zeta) \le \ug_i(y)$
is satisfied; and if equality holds, then
\[
\ug_i(y) = \ug_i(y + s\zeta) = \og_i(y + s\zeta) = \og_i(y + \zeta)
\]
for all $s \in (0, 1)$.
\item \label{item:g.vi} For all $y \in \R^n$ and all $\zeta \in [0, \infty)^n$, the inequalities
$\ug_i(y) \ge \ug_i(y + \zeta)$ and $\og_i(y) \ge \og_i(y + \zeta)$ are satisfied.
\end{enumerate}
\end{lemma}

Now consider the hyperplane $P \subseteq \R^n$ given by
\[
P = \set{y \in \R^n}{y_1 + \dotsb + y_n = 0}
\]
and its unit normal vector
\[
\sigma = \frac{1}{\sqrt{n}} \begin{pmatrix} 1 \\ \vdots \\ 1 \end{pmatrix} \in \R^n.
\]
Let $e_1, \dotsc, e_n$ be the standard basis vectors of $\R^n$
and define
\[
b_i = \sigma - \sqrt{n} e_i \label{def:b}
\]
for $i = 1, \dotsc, n$. Then
\[
|b_i|^2 = n - 1
\]
and
\[
b_i \cdot b_j = -1
\]
for $i \neq j$. Hence $b_1, \dotsc, b_n$ are the corners of a regular $(n - 1)$-simplex in $P$ centred at $0$ with
side length $\sqrt{2n}$. (Indeed the construction is similar to the standard $(n - 1)$-simplex.)
Thus they are the $(n - 1)$-dimensional counterparts to $a_0, \dotsc, a_n$.

Given a function $f \colon P \times \R \to \R$, we write $\tilde{\nabla} f$ for its gradient
with respect to the variable $p \in P$. We want to show the following.

\begin{proposition} \label{prp:projection}
Let $i \in \Z_{n + 1}$. Then there exists a function $f_i \colon P \times \R \to \R$ such that
for almost every $t \in \R$,
\begin{itemize}
\item the function $p \mapsto f_i(p, t)$ belongs to $\BV_\loc^2(P)$ and $\tilde{\nabla} f_i(p, t) \in \{b_1, \dotsc, b_n\}$
for $\Ha^{n - 1}$-almost every $p \in P$; and
\item its graph is $\Gamma_i(t)$, that is, $\Gamma_i(t) = \set{p + f_i(p, t) \sigma}{p \in P}$.
\end{itemize}
\end{proposition}

Before we can prove this result, we need a few lemmas.

\begin{lemma} \label{lem:cones1}
Let $i \in \Z_{n + 1}$.
Suppose that $t \in \R$ and $y_-, y_+ \in \Gamma_i(t)$.
Then
\[
\bigl(y_- + [0, \infty)^n\bigr) \cap \bigl(y_+ - [0, \infty)^n\bigr) \subseteq \Gamma_i(t).
\]
\end{lemma}

\begin{proof}
We first prove that
\[
\bigl(y_- + (0, \infty)^n\bigr) \cap \bigl(y_+ - (0, \infty)^n\bigr) \subseteq \Gamma_i(t).
\]
Let
\[
y \in \bigl(y_- + (0, \infty)^n\bigr) \cap \bigl(y_+ - (0, \infty)^n\bigr).
\]
Define $\zeta_- = y - y_-$ and $\zeta_+ = y_+ - y$. Then $\zeta_-, \zeta_+ \in (0, \infty)^n$.
According to Lemma \ref{lem:g}, this means that
\[
t \ge \ug_i(y_-) \ge \og_i(y_- + \zeta_-) = \og_i(y) \ge \ug_i(y) = \ug_i(y_+ - \zeta_+) \ge \og_i(y_+) \ge t.
\]
Hence $y \in \Gamma_i(t)$. By the semicontinuity of $\ug_i$ and $\og_i$, we also conclude that
\[
\ug_i(y) \le t \le \og_i(t)
\]
for all $y \in \bigl(y_- + [0, \infty)^n\bigr) \cap \bigl(y_+ - [0, \infty)^n\bigr)$.
\end{proof}

\begin{lemma} \label{lem:cones2}
Let $i \in \Z_{n + 1}$. Let $t \in \R$ and $p \in P$. Suppose that
\[
\set{s \in \R}{p + s\sigma \in \Gamma_i(t)} = [s_-, s_+].
\]
Then
\[
\Gamma_i(t) \cap \bigl(p + s_- \sigma - (0, \infty)^n\bigr) = \emptyset
\]
and
\[
\Gamma_i(t) \cap \bigl(p + s_+ \sigma + (0, \infty)^n\bigr) = \emptyset.
\]
\end{lemma}

\begin{proof}
Let $y \in p + s_- \sigma - (0, \infty)^n$. Choose $s < s_-$ such that $y \in p + s\sigma - (0, \infty)^n$ as well.
Then Lemma \ref{lem:g} implies that
\[
\ug_i(y) \ge \og_i(p + s\sigma) \ge \ug_i(p + s\sigma) > t.
\]
Hence $y \not\in \Gamma_i(t)$. The proof of the second statement is similar.
\end{proof}

\begin{lemma} \label{lem:smooth}
There exists a constant $C$ such that the following holds true.
Suppose that $v \colon \R^n \to \R$ is smooth and bounded with $a_j \cdot \nabla v > -1$ for all $j \in \Z_{n + 1}$
and $\sup_{\R^n} |v| \le M$.
Let $i \in \Z_{n + 1}$. Let $\phi \colon P \times \R \to \R$ be the unique function such that
\[
\begin{pmatrix} p + \phi(p, t) \sigma \\ t \end{pmatrix} \in \Phi_i(\graph(v))
\]
for $p \in P$ and $t \in \R$. Then
\begin{equation} \label{eqn:gradient}
|\tilde{\nabla} \phi(p, t)| \le \sqrt{n}
\end{equation}
for all $p \in P$ and $t \in \R$.
Moreover, for any $R > 0$,
\[
\int_{-R}^R \int_{P \cap B_R(0)} |\tilde{\nabla}^2 \phi| \, d\Ha^{n - 1} \, dt \le C \int_{B_{C(M + R)}(0)} |\nabla^2 v| \, dx.
\]
\end{lemma}

Since the proof of this statement is lengthy, we postpone it to the next section. We now prove Proposition \ref{prp:projection}.

\begin{proof}[Proof of Proposition \ref{prp:projection}]
Let $t \in \R$ and $p \in P$. Since $u$ is bounded, the line
\[
\set{t\bfnu_i + \sum_{k = 1}^n (p_k + s\sigma_k) \bfnu_{i + k}}{s \in \R}
\]
must intersect $\graph(u)$. Hence there exists $s \in \R$ with $p + s\sigma \in \Gamma_i(t)$.

If there are $s_-, s_+ \in \R$ with $s_- < s_+$ such that $p + s_- \sigma \in \Gamma_i(t)$ and $p + s_+ \sigma \in \Gamma_i(t)$,
then Lemma \ref{lem:cones1} implies that $\Gamma_i(t)$ has non-empty interior, denoted by
$\mathring{\Gamma}_i(t)$. Because of Lemma \ref{lem:g}.\ref{item:g.v}, we know that $\ug_i(y) = \og_i(y) = t$
for every $y \in \mathring{\Gamma}_i(t)$. Hence for $t_1 \neq t_2$, it follows that
$\mathring{\Gamma}_i(t_1) \cap \mathring{\Gamma}_i(t_2) = \emptyset$. Therefore, there can only be
countably many $t \in \R$ such that $\mathring{\Gamma}_i(t) \neq \emptyset$.
For all other values, we see that $\Gamma_i(t)$ is a graph of a function over $P$. We denote
this function by $f_i(\blank, t)$.

We extend $f_i$ arbitrarily to the remaining values of $t$.

If $t$ is such that $\mathring{\Gamma}_i(t) = \emptyset$, then Lemma \ref{lem:cones2}
shows that for every $y \in \Gamma_i(t)$, the set $\Gamma_i(t)$ is between the cones
$y + (0, \infty)^n$ and $y - (0, \infty)^n$. It follows that $f_i(\blank, t)$ is Lipschitz continuous.

Next we employ an approximation argument in conjunction with Lemma~\ref{lem:smooth}.
Using a standard mollifier, we can find a sequence of smooth, uniformly bounded functions $v_k \colon \R^n \to \R$
such that $v_k \to u$ locally uniformly as $k \to \infty$ and \cite[Proposition 3.7]{Ambrosio-Fusco-Pallara:00}
\[
|D\nabla u|(\Omega) = \lim_{k \to \infty} \int_\Omega |\nabla^2 v_k| \, dx
\]
whenever $\Omega \subseteq \R^n$ is an open, bounded set with $|D\nabla u|(\partial \Omega) = 0$.
It is then easy to modify $v_k$ such that in addition, it satisfies $a_j \cdot \nabla v_k > -1$ in $\R^n$
for every $j \in \Z_{n + 1}$. Hence Lemma \ref{lem:smooth} applies to $v_k$.

From the above convergence, it follows that for any sequence of points $\bfx_k \in \graph(v_k)$, if
$\bfx_k \to \bfx$ as $k \to \infty$, then $\bfx \in \graph(u)$.
If we define $\phi_k$ as in Lemma~\ref{lem:smooth}, then for any fixed $t \in \R$, the functions
$\phi_k(\blank, t)$ are uniformly bounded in $C^{0, 1}(P \cap B_R(0))$ for any $R > 0$. Hence there
is a subsequence that converges locally uniformly. If $t$ is such that $\Gamma_i(t)$ is the graph of $f_i(\blank, t)$,
then it is clear that the limit of any such subsequence must coincide with $f_i(\blank, t)$. Hence in this
case, we have the locally uniform convergence $\phi_k(\blank, t) \to f_i(\blank, t)$ as $k \to \infty$.
The second inequality in Lemma \ref{lem:smooth} implies that
\[
\limsup_{k \to \infty} \int_{-R}^R \int_{P \cap B_R(0)} |\tilde{\nabla}^2 \phi_k| d\Ha^{n - 1} \, dt < \infty
\]
for any $R > 0$. By Fatou's lemma,
\[
\int_{-R}^R \liminf_{k \to \infty} \int_{P \cap B_R(0)} |\tilde{\nabla}^2 \phi_k| d\Ha^{n - 1} \, dt < \infty.
\]
Therefore, for almost every $t \in (-R, R)$, there exists a subsequence $(\phi_{k_\ell}(\blank, t))_{\ell \in \N}$
converging to $f_i(\blank, t)$ locally uniformly and such that
\[
\limsup_{\ell \to \infty} \int_{P \cap B_R(0)} |\tilde{\nabla}^2 \phi_{k_\ell}| d\Ha^{n - 1} < \infty.
\]
We conclude that $f_i(\blank, t) \in \BV_\loc^2(P)$ for almost all $t \in \R$.

We finally need to show that $\tilde{\nabla} f_i(p, t) \in \{b_1, \dotsc, b_n\}$ for almost every $t \in \R$
and $\Ha^{n - 1}$-almost every $p \in P$.

Consider the function $w_i \colon \R^n \to \R$ with
\[
w_i(x) = \frac{u(x) - a_i \cdot x}{\sqrt{n + 1}}, \quad x \in \R^n.
\]
Then for every $t \in \R$,
\[
\begin{split}
\Gamma_i(t) \times \{t\} & = \Phi_i\bigl(\set{\bfx \in \graph(u)}{\bfx \cdot \bfnu_i = t}\bigr) \\
& = \Phi_i\left(\set{\begin{pmatrix} x \\ u(x) \end{pmatrix}}{x \in \R^n \text{ with } w_i(x) = t}\right).
\end{split}
\]
Note further that $\F_i$ coincides up to an $\Ha^n$-null set with $\set{x \in \R^n}{\nabla w_i(x) = 0}$.
Let $\mathcal{Z} \subset \R^n$ denote the set of all points where $u$ is \emph{not} differentiable. By
Rademacher's theorem, this is an $\Ha^n$-null set.
Hence the coarea formula gives
\[
0 = \int_{\F_i \cup \mathcal{Z}} |\nabla w_i| \, dx = \int_{-\infty}^\infty \Ha^{n - 1}\bigl(w_i^{-1}(\{t\}) \cap (\F_i \cup \mathcal{Z})\bigr) \, dt.
\]
In particular, for almost all $t \in \R$,
\[
\Ha^{n - 1}\bigl(w_i^{-1}(\{t\}) \cap (\F_i \cup \mathcal{Z})\bigr) = 0.
\]
As the map $\U$ (defined in the introduction) is Lipschitz
continuous, we conclude that $\U\left(w_i^{-1}(\{t\}) \cap (\F_i \cup \mathcal{Z})\right)$ is an $\Ha^{n - 1}$-null set, too.
Therefore, for $\Ha^{n - 1}$-almost all $y \in \Gamma_i(t)$, the unique point
$x \in \R$ with
\[
\Phi_i(\U(x)) = \begin{pmatrix} y \\ t \end{pmatrix}
\]
belongs to $\R^n \setminus \mathcal{Z}$ and
satisfies $\nabla u(x) \in A \setminus \{a_i\}$.

To put it differently, for almost every $t \in \R$, the following holds true: for $\Ha^{n - 1}$-almost every $p \in P$
the derivative of $u$ exists at the point
\[
\Theta(p, t) = t\nu_i + \sum_{k = 1}^n (p_k + f_i(p, t) \sigma_k) \nu_{i + k}
\]
and belongs to $A \setminus \{a_i\}$.
Furthermore, we know that $f_i(\blank, t)$ is differentiable at $\Ha^{n - 1}$-almost every $p$ by Rademacher's theorem.
At a point $p \in P$ where both statements hold true, we can differentiate the equation
\[
u(\Theta(p, t)) = \frac{t + \sqrt{n} f_i(p, t)}{\sqrt{n + 1}}.
\]
(The right-hand side is the $(n + 1)$-st component of
\[
t\bfnu_i + \sum_{k = 1}^n (p_k + f_i(p, t) \sigma_k) \bfnu_{i + k} = \Phi_i^{-1}\begin{pmatrix} p + f_i(p, t) \sigma \\ t \end{pmatrix}
\]
because $p \in P$ and by the definition of $\sigma$.)
For any $\varpi \in P$, we thus obtain
\[
-\left(\sum_{k = 1}^n \varpi_k a_{i + k} + \varpi \cdot \tilde{\nabla} f_i(p, t) \sum_{k = 1}^n \sigma_k a_{i + k}\right) \cdot \nabla u(\Theta(p, t)) =\sqrt{n} \, \varpi \cdot \tilde{\nabla} f_i(p, t).
\]
If $\nabla u(\Theta(p, t)) = a_j$ for some $j \neq i$, then this simplifies to
\[
- (n + 1) \varpi_{j - i} - \frac{1}{\sqrt{n}} \varpi \cdot \tilde{\nabla} f_i(p, t) = \sqrt{n} \varpi \cdot \tilde{\nabla} f_i(p, t).
\]
Hence
\[
\varpi \cdot \tilde{\nabla} f_i(p, t) = -\sqrt{n} \, \varpi_{j - i} = b_{j - i} \cdot \varpi.
\]
We therefore conclude that $\tilde{\nabla} f_i(p, t) = b_{j - i}$ at such a point.
\end{proof}

\section{Proof of Lemma \ref{lem:smooth}}

In this section we give the postponed proof of Lemma \ref{lem:smooth}. To this end, we first
need another lemma.

\begin{lemma} \label{lem:determinant}
Let $\Lambda$ denote the $(n \times n)$-matrix with columns
\[
\sum_{i \in \Z_{n + 1}} \gamma_{ik} a_i, \quad k = 1, \dotsc, n.
\]
Then
\[
\det(\Lambda) = (-1)^n (n + 1)^{\frac{n - 1}{2}} \det \begin{pmatrix}
\gamma_{01} & \cdots & \gamma_{0n} & 1 \\
\vdots & & \vdots & \vdots \\
\gamma_{n1} & \cdots & \gamma_{nn} & 1
\end{pmatrix}.
\]
\end{lemma}

\begin{proof}
Let $M$ denote the $((n + 1) \times (n + 1))$-matrix with columns
\[
\sum_{i \in \Z_{n + 1}} \gamma_{ik} \bfnu_i, \quad k = 1, \dotsc, n, \quad \text{and} \quad \sum_{i \in \Z_{n + 1}} \bfnu_i.
\]
Then, since $(\bfnu_1, \dotsc, \bfnu_{n + 1})$ is a positively oriented basis of $\R^{n + 1}$, we conclude that
\[
\det(M) = \det\begin{pmatrix}
\gamma_{01} & \cdots & \gamma_{0n} & 1 \\
\vdots & & \vdots & \vdots \\
\gamma_{n1} & \cdots & \gamma_{nn} & 1
\end{pmatrix}.
\]
On the other hand,
\[
M = \frac{1}{\sqrt{n + 1}} \begin{pmatrix}
&&& 0 \\
& -\Lambda & & \vdots \\
&&& 0 \\
m_1 & \cdots & m_n & n + 1
\end{pmatrix},
\]
where $m_k = \sum_{i \in \Z_{n + 1}} \gamma_{ik}$. Hence
\[
\det(M) = (-1)^n (n + 1)^{\frac{1 - n}{2}} \det(\Lambda).
\]
The claim follows immediately.
\end{proof}

\begin{proof}[Proof of Lemma \ref{lem:smooth}]
First we note that by the assumptions on $v$, the intersection of $\graph(v)$ with the hyperplane
$\set{\bfx \in \R^{n + 1}}{\bfx \cdot \bfnu_i = t}$ is a smooth $(n - 1)$-dimensional manifold
for every $t \in \R$. Furthermore, the function $\phi$ is smooth. If we define $\bfXi \colon P \times \R^2 \to \R^{n + 1}$
such that
\[
\bfXi(p, s, t) = t\bfnu_i + \sum_{k = 1}^n (p_k + s \sigma_k) \bfnu_{i + k}
\]
for $p \in P$ and $s, t \in \R$, then $\phi$ is characterised by the condition
that
\[
\bfXi(p, \phi(p, t), t) \in \graph(v)
\]
for all $t \in \R$ and $p \in P$. Hence
\begin{equation} \label{eqn:v-and-phi}
v\bigl(\Xi(p, \phi(p, t), t)\bigr) = \Xi_{n + 1}(p, \phi(p, t), t).
\end{equation}
We now differentiate this equation.

We compute
\[
\dd{\Xi}{t} = \nu_i = - \frac{a_i}{\sqrt{n + 1}}, \quad \dd{\Xi_{n + 1}}{t} = \frac{1}{\sqrt{n + 1}}.
\]
For $\varpi \in P$,
\[
\varpi \cdot \tilde{\nabla} \Xi = -\frac{1}{\sqrt{n + 1}} \sum_{k = 1}^n \varpi_k a_{i + k}, \quad
\varpi \cdot \tilde{\nabla} \Xi_{n + 1} = \frac{1}{\sqrt{n + 1}} \sum_{k = 1}^n \varpi_k = 0.
\]
Finally,
\[
\dd{\Xi}{s} = \sum_{k = 1}^n \sigma_k \nu_{i + k} = - \frac{1}{\sqrt{n^2  + n}} \sum_{k = 1}^n a_{i + k} = \frac{a_i}{\sqrt{n^2 + n}}, \quad
\dd{\Xi_{n + 1}}{s} = \sqrt{\frac{n}{n + 1}}.
\]
We define $\Theta(p, t) = \Xi(p, \phi(p, t), t)$. Differentiating \eqref{eqn:v-and-phi}, we now conclude that
\[
\left(\frac{1}{\sqrt{n}} \dd{\phi}{t}(p, t) - 1\right) a_i \cdot \nabla v(\Theta(p, t)) = \sqrt{n} \dd{\phi}{t}(p, t) + 1
\]
and
\begin{equation} \label{eqn:first-derivative}
\left(\frac{1}{\sqrt{n}} \varpi \cdot \tilde{\nabla} \phi(p, t) a_i - \sum_{k = 1}^n \varpi_k a_{i + k}\right) \cdot \nabla v(\Theta(p, t)) = \sqrt{n} \varpi \cdot \tilde{\nabla} \phi(p, t).
\end{equation}
Hence
\begin{equation} \label{eqn:derivative-of-phi}
\dd{\phi}{t}(p, t) = \sqrt{n} \frac{a_i \cdot \nabla v(\Theta(p, t)) + 1}{a_i \cdot \nabla v(\
\Theta(p, t)) - n}
\end{equation}
and
\[
\varpi \cdot \tilde{\nabla} \phi(p, t) = \sqrt{n} \frac{\sum_{k = 1}^n \varpi_k a_{i + k} \cdot \nabla v(\Theta(p, t))}{a_i \cdot \nabla v(\Theta(p, t)) - n}.
\]

Fix $t \in \R$ and $p \in P$. Since $\nabla v(\Theta(p, t))$ is in the interior of the convex hull of the set $\set{a_j}{j \in \Z_{n + 1}}$,
there exist $\tau_j \in (0, 1)$ for $j \in \Z_{n + 1}$ such that
\[
\sum_{j \in \Z_{n + 1}} \tau_j = 1
\]
and
\[
\nabla v(\Theta(p, t)) = \sum_{j \in \Z_{n + 1}} \tau_j a_j.
\]
Then
\[
a_i \cdot \nabla v(\Theta(p, t)) - n = n \tau_i - \sum_{j \neq i} \tau_j - n = (n + 1) (\tau_i - 1),
\]
while
\[
\sum_{k = 1}^n \varpi_k a_{i + k} \cdot \nabla v(\Theta(p, t)) = \sum_{k = 1}^n \varpi_k \left(n\tau_{i + k} - \sum_{j \neq i + k} \tau_j\right) = (n + 1) \sum_{k = 1}^n \varpi_k \tau_{i + k}.
\]
We further note that
\[
\tau_{i + 1}^2 + \dotsb + \tau_{i + n}^2 \le (\tau_{i + 1} + \dotsb + \tau_{i + n})^2 = (1 - \tau_i)^2.
\]
The Cauchy-Schwarz inequality therefore implies that
\[
\left|\sum_{k = 1}^n \varpi_k a_{i + k} \cdot \nabla v(\Theta(p, t))\right| \le (n + 1) (1 - \tau_i) |\varpi|.
\]
It follows that
\[
|\varpi \cdot \tilde{\nabla} \phi(p, t)| \le \sqrt{n} |\varpi|,
\]
and inequality \eqref{eqn:gradient} is proved.

In order to prove the second statement of Lemma \ref{lem:smooth}, we need to differentiate
\eqref{eqn:first-derivative} again with respect to $p$. We write $\Lambda : M$ for the Frobenius inner product between two matrices
$\Lambda$ and $M$. We also drop the arguments $(p, t)$ in the derivatives of $\phi$ and in $\Theta$.
Then for all $\varpi, \xi \in P$,
\begin{multline*}
\sqrt{\frac{n + 1}{n}} (\xi \otimes \varpi) : \tilde{\nabla}^2 \phi \\
= \frac{\left(\frac{\xi \cdot \tilde{\nabla} \phi}{\sqrt{n}} a_i - \sum_{k = 1}^n \xi_k a_{i +k}\right) \otimes \left(\frac{\varpi \cdot \tilde{\nabla} \phi}{\sqrt{n}} a_i - \sum_{k = 1}^n \varpi_k a_{i +k}\right)}{n - a_i \cdot \nabla v(\Theta)} : \nabla^2 v(\Theta).
\end{multline*}
As we have already seen that $|\tilde{\nabla} \phi| \le \sqrt{n}$, it follows that there is a constant
$C_1 = C_1(n)$ such that
\[
|\tilde{\nabla}^2 \phi| \le \frac{C_1|\nabla^2 v(\Theta)|}{n - a_i \cdot \nabla v(\Theta)}.
\]

Choose an orthonormal basis $(\eta_1, \dotsc, \eta_{n - 1})$ of $P$.
Next we examine the derivative $d\Theta$, and more specifically, its determinant.

Let $\eta_{1k}, \dotsc, \eta_{nk}$ denote the components of $\eta_k$. For $t \in \R$ and $p \in P$, we also define
\[
\eta_{n + 1, k}(p, t) = - \frac{1}{\sqrt{n}} \eta_k \cdot \tilde{\nabla} \phi(p, t), \quad k = 1, \dotsc, n - 1,
\]
and
\[
\eta_{n + 1, n}(p, t) = 1 - \frac{1}{\sqrt{n}} \dd{\phi}{t}(p, t).
\]
Finally, we set $\eta_{\ell n} = 0$ for $\ell = 1, \dotsc, n$. We compute
\[
\eta_k \cdot \tilde{\nabla} \Theta(p, t) = \frac{1}{\sqrt{n + 1}} \left(\frac{1}{\sqrt{n}} \eta_k \cdot \tilde{\nabla} \phi(p, t) a_i - \sum_{\ell = 1}^n \eta_{\ell k} a_{i + \ell}\right)
\]
and
\[
\dd{\Theta}{t}(p, t) = \left(\frac{1}{\sqrt{n}} \dd{\phi}{t}(p, t) - 1\right) \frac{a_i}{\sqrt{n + 1}}.
\]
Hence we can represent $d\Theta$ by the matrix
with columns
\[
-\sum_{\ell = 1}^{n + 1} \frac{\eta_{\ell k} a_{i + \ell}}{\sqrt{n + 1}}, \quad k = 1, \dotsc, n,
\]
with respect to the basis of $P \times \R$ generated by $\eta_1 \dotsc, \eta_{n - 1}$.
Lemma \ref{lem:determinant} now tells us that
\[
\begin{split}
\det(d\Theta) & = \pm \frac{1}{\sqrt{n + 1}} \det\begin{pmatrix}
\eta_{11} & \cdots & \eta_{1n} & 1 \\
\vdots && \vdots & \vdots \\
\eta_{n + 1, 1} & \cdots & \eta_{n + 1, n} & 1
\end{pmatrix} \\
& = \pm \frac{1}{\sqrt{n + 1}} \det \begin{pmatrix}
\eta_{11} & \cdots & \eta_{1, n - 1} & 0 & 1 \\
\vdots && \vdots & \vdots & \vdots \\
\eta_{n1} & \cdots & \eta_{n, n - 1} & 0 & 1 \\
\eta_{n + 1, 1} & \cdots & \eta_{n + 1, n - 1} & \eta_{n + 1, n} & 1
\end{pmatrix} \\
& = \mp \sqrt{\frac{n}{n + 1}} \eta_{n + 1, n} \det\begin{pmatrix}
\eta_{11} & \cdots & \eta_{1, n - 1} & \sigma_1 \\
\vdots && \vdots & \vdots \\
\eta_{n1} & \cdots & \eta_{n, n - 1} & \sigma_n
\end{pmatrix}.
\end{split}
\]
As $(\eta_1, \dotsc, \eta_{n - 1}, \sigma)$ form an orthonormal basis of $\R^n$,
we find that
\[
|\det(d\Theta)| = \sqrt{\frac{n}{n + 1}} |\eta_{n + 1, n}| = \frac{1}{\sqrt{n + 1}} \left|\sqrt{n} - \dd{\phi}{t}\right|.
\]
Recalling \eqref{eqn:derivative-of-phi}, we now obtain
\[
|\det(d\Theta)| = \frac{\sqrt{n^2 + n}}{n - a_i \cdot \nabla v(\Theta)}.
\]
We also note that the map $\Theta$ is injective.
Given $R > 0$, we therefore compute
\begin{equation} \label{eqn:second-derivative}
\begin{split}
\lefteqn{\int_{-R}^R \int_{P \cap B_R(0)} |\tilde{\nabla}^2 \phi| \, d\Ha^{n - 1} \, dt} \qquad\qquad & \\
& \le C_1 \int_{-R}^R \int_{P \cap B_R(0)} \frac{|\nabla^2 v(\Theta)|}{n - a_i \cdot \nabla v(\Theta)} \, d\Ha^{n - 1} \, dt \\
& = \frac{C_1}{\sqrt{n^2 + n}} \int_{-R}^R \int_{P \cap B_R(0)} |\nabla^2 v(\Theta)| |\det(d\Theta)| \, d\Ha^{n - 1} \, dt \\
& = \frac{C_1}{\sqrt{n^2 + n}} \int_{\Theta((P \cap B_R(0)) \times (-R, R))} |\nabla^2 v| \, dx.
\end{split}
\end{equation}

It remains to examine the set $\Theta((P \cap B_R(0)) \times (-R, R))$. Recall that we have the assumption
$\sup_{\R^n} |v| \le M$ in Lemma \ref{lem:smooth}. Thus \eqref{eqn:v-and-phi} implies that
\[
|\Xi_{n + 1}(p, \phi(p, t), t)| \le M.
\]
Since
\[
\Xi_{n + 1}(p, \phi(p, t), t) = \frac{t + \sqrt{n} \phi(p, t)}{\sqrt{n + 1}},
\]
this means that
\[
|\phi(p, t)| \le M \sqrt{\frac{n + 1}{n}} + \frac{R}{\sqrt{n}}
\]
when $t \in (-R, R)$. Hence there exists a constant $C_2 = C_2(n)$ such that
\[
|\Theta(p, t)| \le C_2(M + R)
\]
for all $p \in P \cap B_R(0)$ and all $t \in (-R, R)$.
Thus \eqref{eqn:second-derivative} implies the second inequality of Lemma \ref{lem:smooth}.
\end{proof}

\section{Proof of Theorem \ref{thm:affinely-independent}}

In this section we combine the previous results to prove
the second main theorem. We first consider a function
$u \in \BV_\loc^2(\R^n) \cap L^\infty(\R^n)$ such that $\graph(u)$
is close to the graph of $\lambda_i \wedge \lambda_j$ or
$\lambda_i \vee \lambda_j$ in a cube in $\R^{n + 1}$ with edges parallel to
$\bfnu_1, \dotsc, \bfnu_{n + 1}$. We will give a condition which implies that such a function actually
coincides with $\lambda_i \wedge \lambda_j$ or $\lambda_i \vee \lambda_j$
up to a constant in part of the domain.

For $i, j \in \Z_{n + 1}$ with $i \neq j$ and for $r, R > 0$, we define
\[
Q_{ij}(r, R) = \Biggl\{\sum_{k \in \Z_{n + 1}} c_k \bfnu_k \colon c_i, c_j \in (-r, r) \text{ and } c_k \in (-R, R) \text{ for } k \not\in \{i, j\}\Biggr\}.
\]
Again we consider the map $\U \colon \R^n \to \R^{n + 1}$ with
$\U(x) = (\begin{smallmatrix} x \\ u(x) \end{smallmatrix})$ for $x \in \R^n$. The following is the key statement for the proof
of Theorem \ref{thm:affinely-independent}.

\begin{proposition} \label{prp:affinely-independent}
Let $n \in \N$. For any $\delta > 0$ there exist $\epsilon > 0$ with the following properties.
Let $i, j \in \Z_{n + 1}$ with $i \neq j$. Suppose that $|u(0)| \le \epsilon$
and either
\begin{equation} \label{eqn:wedge}
|u- \lambda_i \wedge \lambda_j| \le \epsilon \quad \text{in $\U^{-1}\bigl(Q_{ij}(1,1)\bigr)$}
\end{equation}
or
\begin{equation} \label{eqn:vee}
|u - \lambda_i \vee \lambda_j| \le \epsilon \quad \text{in $\U^{-1}\bigl(Q_{ij}(1,1)\bigr)$}.
\end{equation}
Then
\begin{equation} \label{eqn:lower-estimate-edge-length}
\Ha^{n - 1}\bigl(\E^* \cap Q_{ij}(\textstyle \frac{1}{4}, 1)\bigr) \ge 2^{n - 1}(1 - \delta).
\end{equation}
If, in addition,
\begin{equation} \label{eqn:upper-estimate-edge-length}
\Ha^{n - 1}\bigl(\E^* \cap Q_{ij}(1, 1)\bigr) \le 2^{n - 1}(1 + \epsilon),
\end{equation}
then there exist $\alpha, \beta \in \R$ such that
\begin{equation} \label{eqn:piecewise-linear}
u = (\lambda_i + \alpha) \wedge (\lambda_j + \beta) \quad \text{in $\textstyle \U^{-1}\bigl(Q_{ij}(\frac{1}{2}, \frac{1}{2})\bigr)$}
\end{equation}
or 
\begin{equation} \label{eqn:piecewise-linear2}
u = (\lambda_i + \alpha) \wedge (\lambda_j + \beta) \quad  \text{in $\textstyle \U^{-1}\bigl(Q_{ij}(\frac{1}{2}, \frac{1}{2})\bigr)$}.
\end{equation}
\end{proposition}

Before we can prove Proposition \ref{prp:affinely-independent}, we need a few more lemmas.
First we need some more information on the functions $f_i$ from
Proposition \ref{prp:projection}. Recall that $f_i(\blank, t) \in \BV_\loc^2(P)$ for almost all $t \in \R$.

Given $i \in \Z_{n + 1}$ and given $t \in \R$ such that $f_i(\blank, t) \in \BV_\loc^2(P)$, let $\mathcal{D}_i'(t) \subseteq P$ denote the approximate jump set
of $\tilde{\nabla} f_i(\blank, t)$. Thus this set is defined analogously to
$\E'$, but for the function $\tilde{\nabla} f_i(\blank, t)$ instead of $u$.
Furthermore, we set
\[
\mathcal{D}_i^\dagger(t) = \set{p + f_i(p, t) \sigma}{p \in \mathcal{D}_i'(t)},
\]
in analogy to $\E^\dagger$.

\begin{lemma} \label{lem:edge-length}
Let $i \in \Z_{n + 1}$. For almost any $t \in \R$,
\[
\mathcal{D}_i^\dagger(t) \times \{t\} \subseteq \Phi_i(\E^* \cup \Null^*).
\]
Hence for any $t_1, t_2 \in R$ and any Borel set $\Omega \subseteq \R^n$,
\[
\int_{t_1}^{t_2} \Ha^{n - 2}(\mathcal{D}_i^\dagger(t) \cap \Omega) \, dt \le \Ha^{n - 1}\bigl(\E^* \cap \Phi_i^{-1}(\Omega \times (t_1, t_2))\bigr).
\]
\end{lemma}

\begin{proof}
Let $p \in P$ and $t \in \R$. Set
\[
\bfx = \Phi_i^{-1}\begin{pmatrix} p + f_i(p, t) \sigma \\ t \end{pmatrix}.
\]
If $\bfx \in \F^*$, then Proposition \ref{prp:convex-independent}
implies that $\graph(u)$ coincides with a hyperplane in a neighbourhood of $\bfx$.
If that hyperplane is perpendicular to $\bfnu_i$, then $p + f_i(p, t) \sigma \in \mathring{\Gamma}_i(t)$
and $t$ belongs to the null set identified in Proposition \ref{prp:projection}. Otherwise, the function $f_i(\blank, t)$
is affine near $p$, and hence $\Phi_i(\bfx)$ cannot belong to $\mathcal{D}_i^\dagger(t) \times \{t\}$.
This implies the first claim.

The second claim is now a consequence of the coarea formula \cite[Theorem 2.93]{Ambrosio-Fusco-Pallara:00}.
\end{proof}

\begin{lemma} \label{lem:g-const}
Let $k \in \{1, \dotsc, n\}$. Suppose that $\underline{s}, \overline{s} \in \R$ with $\underline{s} < \overline{s}$. 
For $z \in \R^{n - 1}$, define $\ell_z(s) = (z_1, \dotsc, z_{k - 1}, s, z_k, \dotsc, z_{n - 1})^T$
for $s \in [\underline{s}, \overline{s}]$, and $L_z = \set{\ell_z(s)}{\underline{s} \le s \le \overline{s}}$.
Fix $i \in \Z_{n + 1}$. Then for $\Ha^{n - 1}$-almost every $z \in \R^{n - 1}$, either
\[
\ug_i(y) = \og_i(y) = \ug_i(y') = \og_i(y')
\]
for all $y, y' \in L_z$, or there exist $\boldsymbol{y} \in L_z \times \R$ such that
\[
\og_i(\ell_z(\overline{s})) \le y_{n + 1} \le \ug_i(\ell_z(\underline{s}))
\]
and $\boldsymbol{y} \in \Phi_i(\E^*)$.
\end{lemma}

\begin{proof}
Consider the projection $\Pi \colon \R^{n + 1} \to \R^n$ given by $\Pi(\boldsymbol{y}) = y$
for $\boldsymbol{y} \in \R^n$. Set $\Psi_i = \Pi \circ \Phi_i$.
Then for $j \in \Z_{n + 1}$ with $j \neq i$ and for $\bfx \in \F_j^*$, it is
clear that $\Jac_{\F^*}\Psi_i(\bfx) = 0$. Hence the area formula gives
$\Ha^n(\Psi_i(\F_j^*)) = 0$. This means that for $\Ha^{n - 1}$-almost every $z \in \R^{n - 1}$,
\begin{equation} \label{eqn:almost-whole-line-not-in-F_j}
\Ha^1(L_z \cap \Psi_i(\F_j^*)) = 0
\end{equation}
for all $j \neq i$.
Furthermore, since $\E^*$ is an $\Ha^{n - 1}$-rectifiable set and $\Ha^{n - 1}(\Null^*) = 0$,
we also know that for $\Ha^{n - 1}$-almost every $z \in \R^{n - 1}$,
\begin{equation} \label{eqn:almost-whole-line-not-in-E}
\Ha^1(L_z \cap \Psi_i(\E^*)) = 0
\end{equation}
and
\begin{equation} \label{eqn:no-N}
L_z \cap \Psi_i(\Null^*) = \emptyset.
\end{equation}

Consider a point $z \in \R^{n - 1}$ such that \eqref{eqn:almost-whole-line-not-in-F_j}, \eqref{eqn:almost-whole-line-not-in-E},
and \eqref{eqn:no-N} hold true.
Recall that by Lemma \ref{lem:g}, a point $\boldsymbol{y} \in \R^{n + 1}$ belongs to
$\Phi_i(\graph(u))$ if, and only if, $\ug_i(y) \le y_{n + 1} \le \og_i(y)$.
Also recall that
\[
\graph(u) = \E^* \cup \Null^* \cup \bigcup_{j \in \Z_{n + 1}} \F_i^*.
\]
From \eqref{eqn:almost-whole-line-not-in-F_j}--\eqref{eqn:no-N} we therefore infer that for $\Ha^1$-almost all $y \in L_z$,
\begin{equation} \label{eqn:in-F_i}
\begin{pmatrix} y \\ t \end{pmatrix} \in \Phi_i(\F_i^*) \quad \text{for all } t \in [\ug_i(y), \og_i(y)].
\end{equation}

Consider $\boldsymbol{y} \in \Phi_i(\F_i^*)$ with $y \in L_z$. Then, setting
$\bfx = \Phi_i^{-1}(\boldsymbol{y})$, we have the locally uniform
convergence $u_{x, \rho} \to \lambda_i$ as $\rho \searrow 0$. Hence for any compact
set $K \subseteq \R^{n + 1}$ and any $\epsilon > 0$ there exists $\rho_0 > 0$ such that
\[
\frac{1}{\rho}(\graph(u) - \bfx) \cap K \subseteq \set{\tilde{\bfx} \in \R^{n + 1}}{\dist(\tilde{\bfx}, \graph(\lambda_i)) < \epsilon/2}
\]
for all $\rho \in (0, \rho_0]$. Recall that $e_1, \dotsc, e_n$ are the standard basis
vectors in $\R^n$. It follows that there exists $r_0 > 0$ such that for all $r \in (0, r_0]$,
\[
\bigl|\ug_i(y \pm r e_k) - \ug_i(y)\bigr| \le r\epsilon \quad \text{and} \quad \bigl|\og_i(y \pm r e_k) - \og_i(y)\bigr| \le r\epsilon
\]
and $|\ug_i(y) - \og_i(y)| \le r\epsilon$.
Thus
\[
\dd{}{y_k}\ug_i(y) = 0 \quad \text{and} \quad \dd{}{y_k}\og_i(y) = 0
\]
and $\ug_i(y) = \og_i(y)$.
Since this is true for $\Ha^1$-almost all $y \in L_z$, Lemma \ref{lem:g}.\ref{item:g.vi} implies that
\begin{equation} \label{eqn:g-inequalities}
\ug_i(\ell_z(\underline{s})) \ge \og_i(y) \ge \ug_i(y) \ge \og_i(\ell_z(\overline{s}))
\end{equation}
for all $y \in L_z$.

If \eqref{eqn:in-F_i} holds for \emph{all} $y \in L_z$, then we immediately conclude that $\ug_i$ and $\og_i$ are
constant and coincide on $L_z$, i.e., we have the first alternative from the statement of the lemma.
If there exists $y \in L_z$ such that \eqref{eqn:in-F_i} does \emph{not} hold true, then
by the above observations, we know that
\[
\begin{pmatrix} y \\ t \end{pmatrix} \not\in \Phi_i(\F_i^*)
\]
holds in fact for all $t \in [\ug_i(y), \og_i(y)]$. Moreover, because \eqref{eqn:in-F_i} still holds true almost
everywhere on $L_z$, there exists a sequence $(\tilde{y}_m)_{m \in \N}$
in $L_z$ such that $y = \lim_{m \to \infty} \tilde{y}_m$ and such that \eqref{eqn:in-F_i} holds for every $\tilde{y}_m$.
We may then choose $\tilde{t}_m \in [\ug_i(\tilde{y}_m), \og_i(\tilde{y}_m)]$. Extracting a subsequence if necessary, we may
assume that $y_{n + 1} = \lim_{m \to \infty} \tilde{t}_m$ exists. Set $\boldsymbol{y} = (\begin{smallmatrix} y \\ y_{n + 1} \end{smallmatrix})$.
Then $\Phi_i^{-1}(\boldsymbol{y})$
belongs to the boundary of $\F_i^*$ relative to $\graph(u)$.

Proposition \ref{prp:convex-independent} implies that $\F_i^*$ is an open set
relative to $\graph(u)$, and its relative boundary is contained in 
$\E^* \cup \Null^*$. Because of \eqref{eqn:no-N}, it follows that
$\Phi_i^{-1}(\boldsymbol{y}) \in \E^*$.
Moreover, \eqref{eqn:g-inequalities} implies that
\[
\og_i(\ell_z(\overline{s})) \le y_{n + 1} \le \ug_i(\ell_z(\underline{s})).
\]
Thus $\boldsymbol{y}$ has the properties from the second alternative in the statement.
\end{proof}

\begin{lemma} \label{lem:g-improved}
Let $i \in \Z_{n + 1}$. Suppose that $G \subseteq \R^n$ is a connected set such that $G \cap \Gamma_i(t) = \emptyset$
for all $t \in (-1, 1)$. Then either $\ug_i(y) \ge 1$ for all $y \in G$ or $\og_i(y) \le -1$ for all $y \in G$.
\end{lemma}

\begin{proof}
Assume that there exists $y_0 \in G$ such that $\ug_i(y_0) < 1$. Since $G \cap \Gamma_i(t) = \emptyset$
for all $t \in (-1, 1)$, this implies that
\[
-1 \ge \og_i(y_0) \ge \ug_i(y_0)
\]
by Lemma \ref{lem:g}.\ref{item:g.iv}.

Given $t \in (-1, 1)$, define
\[
H_t = \set{y \in G}{\og_i(y) \ge t}.
\]
Because $\og_i$ is upper semicontinuous by Lemma \ref{lem:g}, this is a closed set relative to $G$.
Moreover, if $y \in H_t$, it follows that
\[
\og_i(y) \ge \ug_i(y) \ge 1,
\]
because $G \cap \Gamma_i(t') = \emptyset$ for all $t' \in (-1, 1)$. By the lower semicontinuity of
$\ug_i$, this means that there exists $\rho > 0$ such that $\og_i \ge \ug_i \ge t$
in $B_\rho(y)$. Hence $H_t$ is also open relative to $G$.
Since $G$ is connected and $y_0 \not\in H_t$, it follows that $H_t = \emptyset$. This is true for all $t \in (-1, 1)$,
so $\og_i(y) \le -1$ for all $y \in G$.
\end{proof}

We now have everything in place for the proof of Proposition \ref{prp:affinely-independent}.

\begin{proof}[Proof of Proposition \ref{prp:affinely-independent}]
We use induction over $n$. The statement is clear for $n = 1$.
We now assume that $n \ge 2$ and the statement holds true for $n - 1$.

For simplicity, we assume that $i = 1$ and $j = 2$. We also assume that
\eqref{eqn:wedge} holds true; the proof is similar under the assumption
\eqref{eqn:vee}.

Let
\[
\Lambda = \bigl((-\infty, 0] \times \{0\} \times \R^{n - 2}\bigr) \cup \bigl(\{0\} \times (-\infty, 0] \times \R^{n - 2}\bigr).
\]
Then
\[
\Phi_0(\graph(\lambda_1 \wedge \lambda_2)) = \Lambda \times \R.
\]

Let
\[
\epsilon' = \epsilon \sqrt{\frac{n}{n + 1}}.
\]
Under the assumptions of the proposition, the set $\Phi_0(\graph(u)) \cap (-1, 1)^n$
is between $(\Lambda - \epsilon' \sigma) \times \R$ and $(\Lambda + \epsilon' \sigma) \times \R$, i.e.,
\[
\Phi_0(\graph(u)) \cap (-1, 1)^{n + 1} \subseteq \bigcup_{- \epsilon' \le s \le \epsilon'}  (\Lambda + s \sigma) \times \R.
\]

Set $s_0 = \sqrt{\frac{n}{n + 1}} u(0)$. Then $|s_0| \le \epsilon'$ by the assumption that $|u(0)| \le \epsilon$.
Moreover, we compute
\[
\Phi_0\begin{pmatrix} 0 \\ u(0) \end{pmatrix} = \frac{u(0)}{\sqrt{n + 1}} \begin{pmatrix} 1 \\ \vdots \\ 1 \end{pmatrix} = s_0 \begin{pmatrix} \sigma \\ \frac{1}{\sqrt{n}} \end{pmatrix}.
\]
Assuming that $\epsilon < \sqrt{n + 1}$, we infer that
$\og_0(s_0\sigma) > -1$ and $\ug_0(s_0\sigma) < 1$. Using Lemma \ref{lem:g}.\ref{item:g.v} and
Lemma \ref{lem:g-improved}, we conclude that
\[
\ug_0(y) \ge 1 \quad \text{for $y\in (-1, 1)^n \cap \bigcup_{s <- \epsilon'} (\Lambda + s\sigma)$}
\]
and
\[
\og_0(y) \le -1 \quad \text{for $y \in (-1, 1)^n \cap \bigcup_{s > \epsilon'} (\Lambda + s\sigma)$}.
\]

Now consider the function $f_0 \colon P \times \R \to \R$ from
Proposition \ref{prp:projection}. For almost every $t \in (-1, 1)$,
the graph of $f_0(\blank, t)$, which is given by $\Gamma_0(t)$, is between $\Lambda - \epsilon' \sigma$ and
$\Lambda + \epsilon' \sigma$ in the hypercube $(-1, 1)^n$.

Define $\mu_1, \mu_2 \colon P \to \R$
by $\mu_1(p) = b_1 \cdot p$ and $\mu_2(p) = b_2 \cdot p$ for $p \in P$
(where $b_1$ and $b_2$ are the vectors defined on page \pageref{def:b}). Let $F_t \colon P \to \R^n$ be the map with
$F_t(p) = p + f_0(p, t) \sigma$ for $p \in P$. Then
it follows that
\[
|f_0(\blank, t) - \mu_1 \wedge \mu_2| \le \epsilon' \quad \text{in $F_t^{-1}\bigl((-1, 1)^n\bigr)$}.
\]
Moreover, the condition $|f_0(0, t)|\le \epsilon'$ is clearly satisfied.
Hence we may apply the induction hypothesis to the function $f_0(\blank, t)$.
We thereby obtain the inequality
\begin{equation} \label{eqn:lower-estimate-slice}
\Ha^{n - 2}\bigl(\mathcal{D}_0^\dagger(t) \cap \textstyle \bigl((-\frac{1}{4}, \frac{1}{4})^2 \times (-1, 1)^{n - 2}\bigr)\bigr) \ge 2^{n - 2} (1 - \delta)
\end{equation}
for almost all $t \in (-1, 1)$, provided that $\epsilon$ is sufficiently small.
Using Lemma \ref{lem:edge-length}, we therefore obtain inequality \eqref{eqn:lower-estimate-edge-length}.
This proves the first statement of Proposition \ref{prp:affinely-independent}.

In order to prove the second statement, assume now that \eqref{eqn:upper-estimate-edge-length} holds true.
Then
\[
\int_{-1}^1 \Ha^{n - 2}\bigl(\mathcal{D}_0^\dagger(t) \cap (-1, 1)^n\bigr) \, dt \le 2^{n - 1} (1 + \epsilon.)
\]
Recall that we also have inequality \eqref{eqn:lower-estimate-slice}, and we may now assume that $\delta$ is arbitrarily small.
Hence there exist $t_- \in (-1, -\frac{1}{2})$ and
$t_+ \in (\frac{1}{2}, 1)$ such that
\[
\Ha^{n - 2}\bigl(\mathcal{D}_0^*(t_\pm) \cap (-1, 1)^n\bigr) \le 2^{n - 2} (1 + 3\delta + 4\epsilon).
\]
By the induction hypothesis, if $\delta$ and $\epsilon$ are sufficiently small, then
\[
f_0(\blank, t_\pm) = (\mu_1 + \alpha_\pm) \wedge (\mu_2 + \beta_\pm) \quad \text{in $\textstyle F_t^{-1}\bigl((-\frac{1}{2}, \frac{1}{2})^n\bigr)$}
\]
for certain numbers $\alpha_-, \alpha_+, \beta_-, \beta_+ \in \R$.
Therefore, there exist $y_-, y_+ \in \R^2 \times \{0\}^{n - 2}$ such that
\[
\Gamma_0(t_\pm) \cap \textstyle (-\frac{1}{2}, \frac{1}{2})^n = (y_\pm + \Lambda) \cap (-\frac{1}{2}, \frac{1}{2})^n.
\]
Clearly, by the above observations on $\Phi_0(\graph(u))$, this implies that $y_\pm \in B_{\epsilon'}(0)$. We assume that
$\epsilon' \le \frac{1}{4}$.

If $y_- = y_+$, then by Lemma \ref{lem:g},
\[
\Gamma_0(t_\pm) \cap \textstyle (-\frac{1}{2}, \frac{1}{2})^n = (y_+ + \Lambda) \cap (-\frac{1}{2}, \frac{1}{2})^n
\]
for every $t \in (t_-, t_+)$ as well. In this case, we conclude that
\eqref{eqn:piecewise-linear} holds true.
Thus it now suffices to show that $y_- = y_+$.

We argue by contradiction here. Suppose that  $y_- \neq y_+$. We assume
that in fact the first components $y_{1-}$ and $y_{1+}$ are different.
The arguments are similar if $y_{2-} \neq y_{2+}$.

If $y_{1-} \neq y_{1+}$, then for any $z \in (-\frac{1}{2}, -\frac{1}{4}) \times (-\frac{1}{2}, \frac{1}{2})^{n - 2}$, it follows that
\[
\ug_0\begin{pmatrix} y_{1-} \\ z \end{pmatrix} \le t_- \le \og_0 \begin{pmatrix} y_{1-} \\ z \end{pmatrix}
\]
and
\[
\ug_0\begin{pmatrix} y_{1+} \\ z \end{pmatrix} \le t_+ \le \og_0\begin{pmatrix} y_{1+} \\ z \end{pmatrix}.
\]
Since $t_- < t_+$, it is therefore \emph{not} true that $\ug_0$ and $\og_0$
are constant with $\ug_0 = \og_0$ on $[y_{1+}, y_{1-}] \times \{z\}$.
Lemma \ref{lem:g-const} now implies that for $\Ha^{n - 1}$-almost every
$z \in (-\frac{1}{2}, -\frac{1}{4}) \times (-\frac{1}{2}, \frac{1}{2})^{n - 2}$,
the set $[y_{1+}, y_{1-}] \times \{z\} \times [t_-, t_+]$ intersects $\Phi_0(\E^*)$.
It follows that
\[
\Ha^{n - 1}\bigl(\Phi_0(\E^*) \cap (\textstyle (-1, 1) \times (-\frac{1}{2}, -\frac{1}{4}) \times (-1, 1)^{n - 1})\bigr) \ge \frac{1}{4}.
\]
Furthermore, because of \eqref{eqn:lower-estimate-edge-length}, we obtain the estimate
\[
\Ha^{n - 1}\bigl(\E^* \cap Q_{12}(1, 1)\bigr) \ge 2^{n - 1} (1 - \delta) + \frac{1}{4}.
\]
If $\delta + \epsilon < 2^{-n - 1}$, then this contradicts the hypothesis.
\end{proof}

Finally we can prove the second main result with the help of Proposition \ref{prp:convex-independent}
and Proposition \ref{prp:affinely-independent}.

\begin{proof}[Proof of Theorem \ref{thm:affinely-independent}]
Suppose that $A \subseteq \R^n$ is affinely independent. Then $A$ contains at
most $n + 1$ elements. If there are fewer, then we can add additional elements
to $A$ such that it remains affinely independent. Thus we may assume without
loss of generality that the size of $A$ is exactly $n + 1$.

Now suppose that $A = \{\tilde{a}_0, \dotsc, \tilde{a}_n\}$. Consider $M \in \R^{n \times n}$ and $c \in \R^n$ such that $M\tilde{a}_i + c = a_i$
for $i = 0, \dotsc, n$. Then the function $v \colon \R^n \to \R$ with
$v(x) = u(M^T x) + c \cdot x$ has the property that $\nabla v(x) \in \{a_0, \ldots a_n\}$ for almost all $x \in \R^n$. Hence we may assume that $A$ consists
of the vectors $a_0, \dotsc, a_n$. \label{reduction-to-simplex}

Now for the sets $\F$, $\E$, and $\Null$ as defined in Section \ref{sct:faces-edges}, Proposition \ref{prp:convex-independent} implies that $\F \subseteq \mathcal{R}(u)$ with
the same arguments as in the proof of Theorem \ref{thm:convex-independent}.

For $x \in \E$, the functions $u_{x, \rho}$ converge locally uniformly to
$\lambda_i \wedge \lambda_j$ or to $\lambda_i \vee \lambda_j$ as $\rho \searrow 0$
for some $i, j \in \Z_{n + 1}$ with $i \neq j$. Moreover, the approximate tangent
space of $\E^*$ exists at the point $\U(x)$. Clearly this approximate tangent space is $\graph(\lambda_i) \cap \graph(\lambda_j)$. Hence for $\rho$ sufficiently small,
the function $u_{x, \rho}$ satisfies the hypotheses of Proposition \ref{prp:affinely-independent},
including \eqref{eqn:upper-estimate-edge-length}. It follows that $u_{x, \rho}$
satisfies \eqref{eqn:piecewise-linear} or \eqref{eqn:piecewise-linear2}.
In particular, it is regular near $0$, and hence $x \in \mathcal{R}(u)$.

Thus $\mathcal{S}(u) \subseteq \Null$, which is an $\Ha^{n - 1}$-null set.
\end{proof}

\def\cprime{$'$}
\providecommand{\bysame}{\leavevmode\hbox to3em{\hrulefill}\thinspace}
\providecommand{\MR}{\relax\ifhmode\unskip\space\fi MR }
\providecommand{\MRhref}[2]{%
  \href{http://www.ams.org/mathscinet-getitem?mr=#1}{#2}
}
\providecommand{\href}[2]{#2}

\end{document}